\documentclass[11pt]{amsart}

\newtheorem{theorem}[equation]{Theorem}
\newtheorem{cor}[equation]{Corollary}
\newtheorem{prop}[equation]{Proposition}
\newtheorem{lemma}[equation]{Lemma}

\newcommand{\R}{\mathbf{R}}
\newcommand{\Rpq}{\R^{p,q}}
\newcommand{\Rq}{\R^{1,q}}
\newcommand{\Spq}{S^{p,q}}
\newcommand{\Sq}{S^{1,q}}
\newcommand{\Dpq}{\Delta_{\Spq}}
\newcommand{\Dq}{\Delta_{\Sq}}
\newcommand{\x}{\mathbf{x}}
\newcommand{\y}{\mathbf{y}}
\newcommand{\e}{\mathbf{e}}
\newcommand{\C}{\mathbf{C}}

\title[Laplace-Beltrami operator on Hyperboloids]{Eigenfunctions of the Laplace-Beltrami\\operator on hyperboloids}
\author{Amritanshu Prasad}
\address{The Institute of Mathematical Sciences, CIT campus Taramani, Chennai 600113.}
\author{M.~K.~Vemuri}
\address{Chennai Mathematical Institute, Plot H1, SIPCOT IT Park, Padur~PO, Siruseri 603103.}
\begin{document}
\maketitle
\begin{abstract}
  Eigenfunctions of the Laplace-Beltrami operator on a hyperboloid are studied in the spirit of the treatment of the spherical harmonics by Stein and Weiss.
  As a special case, a simple self-contained proof of Laplace's integral for a Legendre function is obtained.
\end{abstract}
In \cite[Chapter~IV, Section~2]{MR0304972}, Stein and Weiss described the spectral decomposition of the Laplace-Beltrami operator on the unit sphere.
Their approach was to identify the eigenfunctions with homogeneous harmonic functions on Euclidean space.

In this article the eigenfunctions of the Laplace-Beltrami operator on a hyperboloid are identified with homogeneous harmonic functions (with respect to a Laplacian of type $(p,q)$) on an open cone.
In the case treated by Stein and Weiss, Liouville's theorem implies that the degree of homogeneity must be a non-negative integer, whereas
here the degree of homogeneity can be any complex number.
This identification is used to compute spherical functions for $O(1,q)$, and consequently Laplace's integral formula for Legendre functions is obtained.
Laplace's integral formula can also be obtained by using the residue theorem \cite[\S15.23]{WW}.
Spherical functions for semisimple Lie groups in general are obtained using different methods (see, e.g., \cite[Chapter~IV]{MR754767}).

Let $n=p+q$.
Let $\Rpq$ denote the space of real $n$-dimensional vectors equipped with the indefinite scalar product of signature $(p,q)$:
\[
\mathbf{x}\cdot \mathbf{y} = {}^t\x Q \y
\]
where $Q$ is the diagonal matrix with $p$ $1$'s followed by $q$ $(-1)$'s along the diagonal. Write $|\x|^2$ for $\x\cdot \x$.
There should be no confusion with the usual positive definite dot product and norm as they are never used in this paper.

Let $\Rpq_+$ denote the subset of $\Rpq$ consisting of those vectors for which $|\x|^2> 0$.
For $\x\in \Rpq_+$, let $|\x|$ denote the positive square root of $|\x|^2$.
Let $O(p,q)$ denote the group consisting of matrices such that ${}^tAQA=Q$.
Denote by $O(p,q)_0$ the connected component of the identity element of $O(p,q)$.
Let $\Spq$ denote the connected component of $(1,0,\ldots,0)$ in the hyperboloid
\[
\{ \x \;:\: |\x|=1,\:x_1>0 \}.
\]
Let $\rho$ be any complex number.
Let $\mathcal{P}_\rho$ denote the space of all functions $f\in C^2(\Rpq_+)$ which are homogeneous of degree $\rho$, i.e., functions such that
\[
f(\lambda \x)=\lambda^\rho f(\x) \mbox{ for all } \x \in \Rpq_+,\: \lambda>0.
\]
Denote by $\Delta$ the differential operator $|\nabla|^2$ (using the indefinite dot product), where $\nabla$ is the gradient operator
\[
\nabla = \left( \frac{\partial}{\partial x_1},\ldots, \frac{\partial}{\partial x_n} \right).
\]
Define
\[
\mathcal{H}_\rho = \{ f\in \mathcal{P}_\rho \;:\: \Delta f = 0 \}.
\]
A function $u\in C^2(\Spq)$ is called a \emph{spherical harmonic\footnote{Perhaps a more apt name would be \emph{hyperboloidal harmonic}.} of degree $\rho$} if $u$ is the restriction to $\Spq$ of a function in $\mathcal{H}_\rho$. Let $H_\rho$ denote the space of spherical harmonics of degree $\rho$:
\[
	 H_\rho = \{ f|_{\Spq} \;: \: f\in \mathcal{H}_\rho \}.
\]
The \emph{Laplace-Beltrami operator} $\Dpq$ on $\Spq$ is defined by
\[
\Dpq u = \Delta \tilde{u}|_{\Spq},
\]
where $\tilde{u}:\Rpq_+\to \C$ is defined by $\tilde{u}(\x)=u(\x/|\x|)$ (the \emph{degree zero homogeneous extension} of $u$). 

Let $\x^\#=Q\x$. The following is easily verified:
\begin{lemma}
  \label{lemma:formulas}
  Let $\x\in \Rpq_+$. Then
  \begin{enumerate}
  \item $\nabla|\x|=\x^\#/|\x|$.
  \item $\nabla|\x|^\rho = \rho |\x|^{\rho -2} \x^\#$.
  \item $|\x^\#|=|\x|$.
  \item $\x^\#\cdot \nabla \tilde{u}(\x) = 0$ for any $u\in C^1(\Spq)$.
  \item $\nabla \cdot \x^\# = n$.
  \end{enumerate}
\end{lemma}
\begin{lemma}
  \label{lemma:eigenvalue}
  If $u\in H_\rho$ then $\Dpq u = -\rho(\rho+n-2)u$.
\end{lemma}
\begin{proof}
  Since $u\in H_\rho$, $|\x|^\rho \tilde{u}(\x) \in \mathcal{H}_\rho$. Therefore (using the formulas in Lemma~\ref{lemma:formulas}),
  \begin{eqnarray*}
    0 &=&\Delta(|\x|^\rho \tilde{u}(\x))\\
    &=&\nabla\cdot(\nabla (|\x|^\rho \tilde{u}(\x)))\\
    &=&\nabla\cdot(\rho|\x|^{\rho-2}\x^\# \tilde{u}(\x)+|\x|^\rho \nabla \tilde{u}(\x))\\
    &=&(\nabla(\rho|\x|^{\rho-2}\tilde{u}(\x))\cdot \x^\# + \rho|\x|^{\rho-2}\tilde{u}(\x)(\nabla \cdot \x^\#) + |\x|^\rho \Delta \tilde{u}(\x)\\
    &=&\rho(\rho-2)|\x|^{\rho-4}\tilde{u}(\x)|\x^\#|^2 + \rho|\x|^{\rho-2} \nabla \tilde{u}(\x)\cdot \x^\# + n\rho|\x|^{\rho-2}\tilde{u}(\x)+\Delta \tilde{u}(\x)\\
    &=&\rho(\rho-2)|\x|^{\rho-4}\tilde{u}(\x)|\x|^2 + n\rho|\x|^{\rho-2}\tilde{u}(\x)+\Delta \tilde{u}(\x).
  \end{eqnarray*}
  Setting $|\x|=1$ in the result of the above calculation yields
  \[
  0=\rho(\rho-2+n)\tilde{u}(\x)+\Delta \tilde{u}(\x),
  \]
from which the lemma follows.
\end{proof}
The following proposition gives a construction of spherical harmonics when $p=1$:
\begin{prop}
  \label{prop:harmonic}
 Suppose $\mathbf{c}\in \mathbf{R}^{1,q}_+$ is an isotropic vector, (meaning that $|\mathbf{c}|^2=0$) such that $c_1>0$. Then $\mathbf{c}\cdot \x>0$ for all $\x\in \Sq$.
 Let $f(\x)=(\mathbf{c}\cdot \x)^\rho$.
 Then $f\in \mathcal{H}_\rho$.
\end{prop}
\begin{proof}
  The set of points where $\mathbf{c}\cdot \x=0$ form a hyperplane tangential to the cone $|\mathbf{c}|^2=0$.
  For fixed $\mathbf{x}$, the sign of $\mathbf{c}\cdot \x$ can change only when $\mathbf{c}$ crosses this hyperplane.
  However, the entire half-cone
  \begin{equation*}
    \{\mathbf{c}\;:\:|\mathbf{c}|^2=0, \:c_1>0\}
  \end{equation*}
  lies on one side of the hyperplane, because the cone is quadratic.
  Therefore, for each $\x\in S^{1,q}$, it suffices to verify that $\mathbf{c}\cdot \x>0$ for $\mathbf{c}=(1,1,0,\ldots,0)$.
  In this case, $\mathbf{c}\cdot \x=x_1-x_2$, which is positive since $x_1>0$ and $x_1^2-x_2^2-\cdots-x_n^2=1$, so that $x_1>|x_i|$ for each $i>1$.

  If $g\in C^2(\Rpq_+)$ and $\phi\in C^2(\R)$, then 
  \begin{equation*}
    \Delta (\phi\circ g)(\x) = \phi''(g(\x))|\nabla g(\x)|^2 + \phi'(g(\x))\Delta g(\x).
  \end{equation*}
  Let $g(\x)=\mathbf{c}\cdot \x$, then $\nabla g(\x)=\mathbf{c}$, so that $|\nabla g(\x)|^2=0$. Since $g$ is linear, $\Delta g(\x)=0$. 
  Therefore $\Delta f(\x)=0$.
\end{proof}
Let $\e= (1,0,\ldots,0)$.
Then $K=\mathrm{Stab}_{O(1,q)_0}(\e)$ is isomorphic to $SO(q)$ and is a maximal compact subgroup of $O(1,q)_0$.
The action of $O(1,q)_0$ on $\Sq$ is transitive, and the $K$-invariant spherical harmonics on $\Sq$ are precisely the $K$-invariant spherical functions for $O(1,q)_0$.
It follows from Proposition~\ref{prop:harmonic} that
\begin{prop}
  \label{prop:invariant_harmonic}
  Let $\mathbf{c}$ be any isotropic vector in $\Rq$. Then
  \begin{equation*}
    \int_K (k\mathbf{c}\cdot \x)^\rho dk
  \end{equation*}
  is a $K$-invariant spherical harmonic of degree $\rho$ on $\Sq$.
\end{prop}
Since $K$ acts transitively on the slices of $\Sq$ by the hyperplanes on which the first coordinate $x_1$ is constant, the value of a $K$-invariant spherical harmonic is simply a function of $x_1$, which will be denoted by $P(x_1)$.
A $K$-invariant spherical harmonic may be viewed as a solution to an ordinary differential equation in $x_1$:
\begin{theorem}
  \label{theorem:ODE}
  Suppose that $P_\rho(x_1)$ is the value of a $K$-invariant spherical harmonic which is homogeneous of degree $\rho$.
  Then $P_\rho$ is a solution to the differential equation
  \begin{equation}
    \label{eq:ODE}
    (1-x_1^2)P''_\rho(x_1)+(1-n)x_1P'(x_1)+\rho(\rho-2+n)P(x_1)=0.
  \end{equation}
\end{theorem}
\begin{proof}
  For any $f\in C^2(\Sq)$
  we have
  \begin{eqnarray*}
    \nabla f(x_1/|\x|) & = & \nabla(x_1/|\x|)f'(x_1/|\x|)\\
    &=& \frac{(\nabla x_1)|\x|-x_1\nabla|\x|}{|\x|^2}f'(x_1/|\x|)\\
    &=& \frac{\e|\x|-(x_1/|\x|)\x^\#}{|\x|^2}f'(x_1/|\x|)\\
    &=& u\mathbf{v},
  \end{eqnarray*}
  where $u=|\x|^{-3}f'(x_1/|\x|)$ and $\mathbf{v}=\e|\x|^2-x_1\x^\#$.
  Since $\Delta=|\nabla|^2$, 
  \begin{equation}
    \label{eq:laplacian}
    \Delta f(x_1/|\x|)=(\nabla u) \cdot \mathbf{v}+ u \nabla\cdot \mathbf{v}.
  \end{equation}
  Now,
  \begin{eqnarray*}
    \nabla u & = & -3|\x|^{-5}\x^\#f'(x_1/|\x|)+|\x|^{-3}\nabla(x_1/|\x|)f''(x_1/|\x|)\\
    &=& -3|\x|^{-5}\x^\#f'(x_1/|\x|)+|\x|^{-6}(\e|\x|^2-x_1\x^\#)f''(x_1/|\x|).
  \end{eqnarray*}
  and
  \begin{equation*}
    \nabla \cdot \mathbf{v} = \e\cdot \nabla |\x|^2-(\e\cdot \x^\#+nx_1) = (1-n)x_1.
  \end{equation*}
  Suppose there exists a function $P$ such that $P(x_1)=f(\x)$ for each $\x$ such that $|\x|=1$.
  Substituting the above values of $\nabla u$ and $\nabla\cdot \mathbf{v}$ in (\ref{eq:laplacian}) and then setting $|\x|=1$ we have, 
  \begin{multline*}
    \Dq f|_{\Sq}(\x)=(-3\x^\#P'(x_1)+(\e-x_1\x^\#)P''(x_1))\cdot(\e-x_1\x^\#)\\
    + (1-n)x_1P'(x_1).
  \end{multline*}
  When $|\x|=1$, $|\e-x_1\x^\#|^2=(1-x_1^2)$ and $(\e-x_1\x^\#)\cdot \x^\#=0$ so that the above equality simplifies to
  \begin{equation*}
    \Dq f|_{\Sq}(\x)= (1-x_1^2)P''(x_1)+(1-n)x_1P'(x_1).
  \end{equation*}
  Combining this with Lemma~\ref{lemma:eigenvalue} gives (\ref{eq:ODE}).
\end{proof}
\begin{cor}
  \label{cor:uniqueness}
  For $n\geq 3$, there is (up to scaling) a unique $K$-invariant spherical function of degree $\rho$ given by
  \begin{equation*}
    \int_K(k\mathbf{c}\cdot \x)^\rho dk,
  \end{equation*}
where $\mathbf{c}$ is any non-zero isotropic vector in $\Rq$.
\end{cor}
\begin{proof}
  The ordinary differential equation (\ref{eq:ODE}) is linear of degree 2 with a regular singular point at $x_1=1$. 
  The indicial equation at this point is
  \begin{equation*}
    m(m+(n-1)/2-1)=0.
  \end{equation*}
  Therefore, it has (up to scaling) at most one solution defined on $[1,\infty)$.
  This solution is known by Proposition~\ref{prop:invariant_harmonic}.
\end{proof}
The classical integral formula due to Laplace for Legendre functions is readily derived from the preceding analysis:
\begin{cor}
  \label{cor:Laplace_integral}
  Every solution of the ordinary differential equation
  \begin{equation*}
    (1-x^2)P''(x)-2xP'(x)+\rho(\rho+1)P(x)=0
  \end{equation*}
  that is defined on $[1,\infty)$ is a scalar multiple of 
  \begin{equation*}
    P_\rho(x)=\frac{1}{2\pi}\int_0^{2\pi} (x+\sqrt{x^2-1}\,\cos\theta)^\rho d\theta.
  \end{equation*}
\end{cor}
\begin{proof}
  Evaluate the formula from Corollary~\ref{cor:uniqueness} taking $q=2$, $\mathbf{c}=(1,0,-1)$ and $\mathbf{x}=(x,0,\sqrt{x^2-1})$.
\end{proof}
\bibliographystyle{amsalpha}
\bibliography{refs}
\end{document}